\documentclass[12pt]{article}
\makeatletter \oddsidemargin  -.1in \evensidemargin -.1in
\newcommand{\singlespacing}{\let\CS=\@currsize\renewcommand{\baselinestretch}{1}\tiny\CS}
\newcommand{\oneandahalfspacing}{\let\CS=\@currsize\renewcommand{\baselinestretch}{1.25}\tiny\CS}
\newcommand{\doublespacing}{\let\CS=\@currsize\renewcommand{\baselinestretch}{1.35}\tiny\CS}

\newtheorem{rule-def}[theorem]{Rule}

\RequirePackage[dvips]{graphicx} \textheight 22.5cm
\usepackage{latexsym,epsfig,enumerate,amsmath,amsfonts,amssymb,amsbsy,amsopn,mathrsfs}
\usepackage{subfigure}
\usepackage{graphicx}

\begin{document}
\title{\bf A Study on Erd\H{o}s-Straus conjecture on Diophantine
  equation $\frac{4}{n}=\frac{1}{x}+\frac{1}{y}+\frac{1}{z}$} \author{\small S. Maiti$^{1,2}$
  \thanks{Corresponding author, Email address: {\it
      maiti0000000somnath@gmail.com/somnathm@lnmiit.ac.in
      (S. Maiti)}}\\\it $^{1}$ Department of Mathematics, The LNM
  Institute of Information Technology\\\it Jaipur 302031, India\\ \it
  $^{2}$Department of Mathematical Sciences, Indian Institute of
  Technology (BHU),\\ Varanasi-221005, India} \date{}
\maketitle \noindent \doublespacing
\vspace{-0.5cm}
\begin{abstract}
The Erd\H{o}s-Straus conjecture is a renowned problem which describes
that for every natural number $n~(\ge 2)$, $\frac{4}{n}$ can be
represented as the sum of three unit fractions. The main purpose of
this study is to show that the Erd\H{o}s-Straus conjecture is
true. The study also re-demonstrates Mordell theorem which states that
$\frac{4}{n}$ has a expression as the sum of three unit fractions for
every number $n$ except possibly for those primes of the form
$n\equiv r$ (mod 780) with $r=1^2,11^2,13^2,17^2,19^2,23^2$. For
$l,r,a\in\mathbb{N}$;
$\frac{4}{24l+1}-\frac{1}{6l+r}=\frac{4r-1}{(6l+r)(24l+1)}$ with $1\le
r\le 12l$, if at least one of the sums in right side of the
expression, say, $a+(4r-a-1),~1\le a\le 2r-1$ for at least one of the
possible value of $r$ such that $a,(4r-a-1)$ divide $(6l+r)(24l+1)$;
then the conjecture is valid for the corresponding $l$. However, in
this way the conjecture can not be proved only twelve values of $l$
for $l$ up to $l=10^5$.

\end{abstract}

\section{Introduction}
The famous Erd\H{o}s-Straus conjecture in number theory, formulated by Paul
Erd\H{o}s and Ernst G. Strauss \cite{Elsholtz,Kotsireas} in 1948,
states that for every natural number $n~(\ge 2)$, $\exists$ natural
numbers $x,y,z$ such that $\frac{4}{n}$ can be expressed as
\begin{equation}
\frac{4}{n}=\frac{1}{x}+\frac{1}{y}+\frac{1}{z}.
\label{Erdos-Straus_problem}
\end{equation}
Many researchers, not only in Number Theory but also in different areas
of Mathematics, gave attention to this conjecture such as L. Bernstein
\cite{Bernstein}, M. B. Crawford \cite{Crawford}, M. Di Giovanni,
S. Gallipoli, M. Gionfriddo \cite{Giovanni}, C. Elsholtz and T.  Tao
\cite{Elsholtz2}, J. Ghanouchi \cite{Ghanouchi,Ghanouchi2},
L. J. Mordell \cite{Mordell}, D. J. Negash \cite{Negash},
R. Obl\'{a}th \cite{Oblath}, L. A. Rosati \cite{Rosati}, J. W. Sander
\cite{Sander}, R. C. Vaughan \cite{Vaughan}, K. Yamamoto
\cite{Yamamoto}, J. W. Porras Ferreira \cite{Ferreira1,Ferreira2} etc. The
validity of the conjecture for all $n\le 10^{14}$ and $n\le 10^{17}$
was reported by Swett \cite{Swett} and Salez \cite{Salez}
respectively.

If $m$ and $n$ are relatively prime integers, then Schinzel
\cite{Schinzel} established that
$\frac{4}{mt+n}=\frac{1}{x(t)}+\frac{1}{y(t)}+\frac{1}{z(t)}$ having
$x(t)$, $y(t)$ and $z(t)$ as integer polynomials in $t$ together with
positive leading coefficients and non quadratic residue $n$ (mod
$m$). Mordell \cite{Mordell} demonstrated that the validity of the
Erd\H{o}s-Straus conjecture for all $n$ except possible cases where
$n$ is congruent to $1^2$, $11^2$, $13^2$, $17^2$, $19^2$ or $23^2$
(mod 840).

The conjecture can be proved if it is derived for the all prime
numbers to any of the following cases (i) $n=p=2m+1$ (ii) $n=p=4m+1$
(iii) $n=p=8m+1$ (iv) $n=p=24m+1$ where $m\in\mathbb{N}$. For
$l,r,a\in\mathbb{N}$;
$\frac{4}{24l+1}-\frac{1}{6l+r}=\frac{4r-1}{(6l+r)(24l+1)}$ with $1\le
r\le 12l$. If at least one of the sums in right side of the
expression, say, $a+(4r-a-1),~1\le a\le 2r-1$ for at least one of the
possible value of $r$ such that $a,(4r-a-1)$ divide $(6l+r)(24l+1)$;
then $\frac{4}{24l+1}$ has the expression in the form of the equation
(\ref{Erdos-Straus_problem}). However, in this way the expression in
the form of equation (\ref{Erdos-Straus_problem}) can not be proved
for all values of $l$ although it can be established for most of the
values of $l$. The computations of
$\frac{4}{24l+1}-\frac{1}{6l+r}=\frac{4r-1}{(6l+r)(24l+1)}$ for $l$ up
to $l=10^5$ has been carried out which shows that the expression in
the form of the equation (\ref{Erdos-Straus_problem}) with $n=24l+1$
can not be proved only twelve values of $l$. Finally, it has been
shown (by other way) in Section \ref{Erdos-Straus_3rd_theorem} that
the conjecture is true.

\section{Results and Discussion}
\subsection{Demonstration of solutions of Erdos-Straus Problem}
If we want to express any $n\in \mathbb{N}$ in the form of
(\ref{Erdos-Straus_problem}), then it is equivalent to find the
equation (\ref{Erdos-Straus_problem}) for all primes $n=p$. We know
that (i) if $n=2m$, then
$\frac{4}{n}=\frac{4}{2m}=\frac{1}{2m}+\frac{1}{2m}+\frac{1}{m}$ (ii)
if $n=3m$, then
$\frac{4}{n}=\frac{4}{3m}=\frac{1}{2m}+\frac{1}{2m}+\frac{1}{3m}=\frac{1}{m}+\frac{1}{3(m+1)}+\frac{1}{3m(m+1)}=\frac{1}{3m}+\frac{1}{m+1}+\frac{1}{m(m+1)}$
(iii) if $n=3m+2$, then
$\frac{4}{n}=\frac{4}{3m+2}=\frac{1}{3m+2}+\frac{1}{m+1}+\frac{1}{(m+1)(3m+2)}$
(iv) if $n=4m+3$, then
$\frac{4}{n}=\frac{4}{4m+3}=\frac{1}{m+1}+\frac{1}{2(4m+3)(m+1)}+\frac{1}{2(4m+3)(m+1)}$
(v) if $n=4m$ and $n=4m+2$, then these reduce to case (i) where $m\in
\mathbb{N}$. Thus we have to prove the equation
(\ref{Erdos-Straus_problem}) for all $n=4m+1$.

From the equation (\ref{Erdos-Straus_problem}), we get
$\frac{4}{n}=\frac{1}{x}+\frac{1}{y}+\frac{1}{z}>\frac{1}{x},\frac{1}{y},\frac{1}{z}$
i.e. $\frac{1}{x},\frac{1}{y},\frac{1}{z}\le
\frac{1}{[\frac{n}{4}]+1}$ as $\frac{4}{n}-\frac{1}{[\frac{n}{4}]}\le
0$. If $x\le y\le z$, then $\frac{4}{n}\le\frac{3}{x}$
i.e. $\frac{1}{[\frac{3n}{4}]}\le\frac{1}{x}$. Thus
$\frac{1}{[\frac{3n}{4}]}\le\frac{1}{x}\le\frac{1}{[\frac{n}{4}]+1}$ and
$\frac{1}{z}\le\frac{1}{y}\le\frac{1}{x}$.

For $\frac{4}{n}-\frac{1}{x}=\frac{4x-n}{nx}$ with $x\le y\le z$ and
$\frac{1}{[\frac{3n}{4}]}\le\frac{1}{x}\le\frac{1}{[\frac{n}{4}]+1}$,
if $4x-n=p_{i_1}^{\alpha_1}p_{i_2}^{\alpha_2}\cdots
p_{i_r}^{\alpha_r}+p_{j_1}^{\beta_1}p_{j_2}^{\beta_2}\cdots
p_{j_s}^{\beta_s}$ with $nx=p_1^{\gamma_1}p_2^{\gamma_2}\cdots
p_m^{\gamma_m},~\{p_{i_1},p_{i_2},\cdots,p_{i_r};p_{j_1},p_{j_2},\cdots,p_{j_s}\}\subset
\{p_1,p_2,\cdots,p_{m}\}$; $p_{i_1}^{\alpha_1}p_{i_2}^{\alpha_2}\cdots
p_{i_r}^{\alpha_r},p_{j_1}^{\beta_1}p_{j_2}^{\beta_2}\cdots
p_{j_s}^{\beta_s}$ divide $nx$; then $n$ has a expression as in the form of
(\ref{Erdos-Straus_problem}).

Again if $4x-n=y_1+z_1;~y_1,z_1\in\mathbb{N}$ with $(y_1,z_1)=1$, then
$\frac{4}{n}-\frac{1}{x}=\frac{4x-n}{nx}=\frac{y_1+z_1}{nx}=\frac{1}{\frac{nx}{y_1}}+\frac{1}{\frac{nx}{z_1}}$. If
$n$ satisfies the equation (\ref{Erdos-Straus_problem}), then
$nx=k_1y_1,~k_1\in\mathbb{N}$ and $nx=k_2z_1,~k_2\in\mathbb{N}$
i.e. $k_1y_1=k_2z_1$. Or, $k_2=gy_1,~g\in\mathbb{N}$ as
$(y_1,z_1)=1$. Then $k_1=gz_1$ and $nx=gy_1z_1$. Thus
$\frac{4}{n}=\frac{1}{x}+\frac{1}{gz_1}+\frac{1}{gy_1}=\frac{1}{x}+\frac{1}{y}+\frac{1}{z}$
with $g=\frac{nx}{y_1z_1}$, $y=gy_1$ and $z=gz_1$.

If
$4x-n=dy_1+dz_1;~d,y_1,z_1\in\mathbb{N}$ with $(y_1,z_1)=1$, then
$\frac{4}{n}-\frac{1}{x}=\frac{4x-n}{nx}=\frac{dy_1+dz_1}{nx}=\frac{1}{\frac{nx}{dy_1}}+\frac{1}{\frac{nx}{dz_1}}$. If
$n$ satisfies the equation (\ref{Erdos-Straus_problem}), then
$nx=dk_1y_1,~k_1\in\mathbb{N}$ and $nx=dk_2z_1,~k_2\in\mathbb{N}$
i.e. $k_1y_1=k_2z_1$. Or, $k_2=gy_1,~g\in\mathbb{N}$. Then $k_1=gz_1$
and $nx=gdy_1z_1$. Thus
$\frac{4}{n}=\frac{1}{x}+\frac{1}{gz_1}+\frac{1}{gy_1}=\frac{1}{x}+\frac{1}{y}+\frac{1}{z}$
with $g=\frac{nx}{dy_1z_1}$, $y=gy_1$ and $z=gz_1$.

For $\frac{4}{4m+1}$, $[\frac{4m+1}{4}]=m$ and
$[\frac{3(4m+1)}{4}]=3m$, $\frac{1}{3m}\le \frac{1}{x}\le
\frac{1}{m+1}$. Then for the equation (\ref{Erdos-Straus_problem}) with $n=4m+1$,
the possible cases of $\frac{1}{x}$ will be only
$\frac{1}{x}=\frac{1}{m+1},~\frac{1}{m+2},~\cdots,~\frac{1}{3m}$. Thus
for
\begin{equation}
\frac{4}{4m+1}-\frac{1}{m+r}=\frac{4r-1}{(m+r)(4m+1)}=\frac{\{1+(4r-2)\},\{2+(4r-3)\},\cdots,\{(2r-1)+2r\}}{(m+r)(4m+1)}
\label{Erdos-Straus_problem_x_and_sum}
\end{equation}
with $r\in\mathbb{N}$ and $1\le r\le 2m$; if at least one of the sums
in right side of the equation, say, $a+(4r-a-1),~1\le a\le 2r-1$ for
at least one of the possible value of $r$ such that $a,(4r-a-1)$
divide $(m+r)(4m+1)$; then $\frac{4}{4m+1}$ satisfies the equation
(\ref{Erdos-Straus_problem}). A question is arising naturally that
whether the equation (\ref{Erdos-Straus_problem}) can be proved with
the help of the equation (\ref{Erdos-Straus_problem_x_and_sum})?
The answer of this question will be given later.

If $m=2k-1$, then $4m+1=8k-3$,
$[\frac{8k-3}{4}]=2k-1,~[\frac{3(8k-3)}{4}]=6k-3,~\frac{1}{6k-3}\le\frac{1}{x}\le
\frac{1}{2k}$. So,
$\frac{4}{8k-3}-\frac{1}{2k}=\frac{3}{2k(8k-3)}=\frac{1+2}{2k(8k-3)}=\frac{1}{2k(8k-3)}+\frac{1}{k(8k-3)}$
i.e. $\frac{4}{8k-3}=\frac{1}{2k}+\frac{1}{2k(8k-3)}+\frac{1}{k(8k-3)}$.

If $m=2k$, then $4m+1=8k+1$,
$[\frac{8k+1}{4}]=2k,~[\frac{3(8k+1)}{4}]=6k,~\frac{1}{6k}\le\frac{1}{x}\le
\frac{1}{2k+1}$. So,
$\frac{4}{8k+1}-\frac{1}{2k+1}=\frac{3}{(2k+1)(8k+1)}$. 

Again if $k=3l-2$, then $8k+1=24l-15$,
$[\frac{24l-15}{4}]=6l-4,~[\frac{3(24l-15)}{4}]=18l-12,~\frac{1}{18l-12}\le\frac{1}{x}\le
\frac{1}{6l-3}$. So,
$\frac{4}{24l-15}-\frac{1}{6l-3}=\frac{3}{(6l-3)(24r-15)}=\frac{1}{(2l-1)(24l-15)}=\frac{1}{2(2l-1)(24l-15)}+\frac{1}{2(2l-1)(24l-15)}$
i.e. $\frac{4}{24l-15}=\frac{1}{6l-3}+\frac{1}{2(2l-1)(24l-7)}+\frac{1}{2(2l-1)(24l-15)}$.

If $k=3l-1$, then $8k+1=24l-7$,
$[\frac{24l-7}{4}]=6l-2,~[\frac{3(24l-7)}{4}]=18l-6,~\frac{1}{18l-6}\le\frac{1}{x}\le
\frac{1}{6l-1}$. So,
$\frac{4}{24l-7}-\frac{1}{6l-1}=\frac{3}{(6l-1)(24l-7)}=\frac{6l}{2l(6l-1)(24l-7)}=\frac{1+(6l-1)}{2l(6l-1)(24l-7)}=\frac{1}{2l(6l-1)(24l-7)}+\frac{1}{2l(24l-7)}$
i.e. $\frac{4}{24l-7}=\frac{1}{6l-1}+\frac{1}{2l(6l-1)(24l-7)}+\frac{1}{2l(24l-7)}$.

If $k=3l$, then $8k+1=24l+1$,
$[\frac{24l+1}{4}]=6l,~[\frac{3(24l+1)}{4}]=18l$. Thus we have to
prove the equation (\ref{Erdos-Straus_problem}) for all $n=24l+1$
with $\frac{1}{18l}\le\frac{1}{x}\le \frac{1}{6l+1}$. Using the
expression of the equation (\ref{Erdos-Straus_problem_x_and_sum}), we get
\begin{equation}
\frac{4}{24l+1}-\frac{1}{6l+r}=\frac{4r-1}{(6l+r)(24l+1)}=\frac{\{1+(4r-2)\},\{2+(4r-3)\},\cdots,\{(2r-1)+2r\}}{(6l+r)(24l+1)}
\label{Erdos-Straus_problem_x_and_sum_24l_sum_1}
\end{equation}
with $r\in\mathbb{N}$ and $1\le r\le 12l$. If at least one of the sums
in right side of the expression, say, $a+(4r-a-1),~1\le a\le 2r-1$ for at least one of
the possible value of $r$ such that $a,(4r-a-1)$ divide
$(6l+r)(24l+1)$; then $\frac{4}{24l+1}$ has the expression in the form
of the equation (\ref{Erdos-Straus_problem}). Thus the same question is arising
naturally that whether the equation (\ref{Erdos-Straus_problem}) can
be proved with the help of the equation
(\ref{Erdos-Straus_problem_x_and_sum_24l_sum_1})?

The answer of this question is that the expression in the form of
equation (\ref{Erdos-Straus_problem}) can not be proved with the help
of the equation (\ref{Erdos-Straus_problem_x_and_sum_24l_sum_1}) for
all values of $l\in\mathbb{N}$. However, the equation
(\ref{Erdos-Straus_problem}) can be proved with the help of the
equation (\ref{Erdos-Straus_problem_x_and_sum_24l_sum_1}) for most of
the values of $l\in\mathbb{N}$ and it can not be proved for very less
number of values of $l\in\mathbb{N}$. To find the answer, the
computations of $l$ up to $l=10^5$ has been carried out which shows
that the expression in the form of the equation
(\ref{Erdos-Straus_problem}) with $n=24l+1$ can not be proved with the
help of the equation (\ref{Erdos-Straus_problem_x_and_sum_24l_sum_1})
only if $\frac{4}{n}=\frac{4}{409}$ (i.e. for $l=17$),
$\frac{4}{n}=\frac{4}{577}$ (i.e. for $l=24$),
$\frac{4}{n}=\frac{4}{5569}$ (i.e. for $l=232$),
$\frac{4}{n}=\frac{4}{9601}$ (i.e. for $l=400$),
$\frac{4}{n}=\frac{4}{23929}$ (i.e. for $l=997$),
$\frac{4}{n}=\frac{4}{83449}$ (i.e. for $l=3477$),
$\frac{4}{n}=\frac{4}{102001}$ (i.e. for $l=4250$),
$\frac{4}{n}=\frac{4}{329617}$ (i.e. for $l=13734$),
$\frac{4}{n}=\frac{4}{712321}$ (i.e. for $l=29680$),
$\frac{4}{n}=\frac{4}{1134241}$ (i.e. for $l=47260$),
$\frac{4}{n}=\frac{4}{1724209}$ (i.e. for $l=71842$),
$\frac{4}{n}=\frac{4}{1726201}$ (i.e. for $l=71925$) and the
computations are taking several days to give result in Mathematica in
my system with 8 GB RAM and 1 TB har disk if the equation
(\ref{Erdos-Straus_problem}) with $n=24l+1$ can not be proved with the
help of the equation (\ref{Erdos-Straus_problem_x_and_sum_24l_sum_1})
for $l\ge 29680$.

Thus for exceptional cases, we have to find the expression of the
equation (\ref{Erdos-Straus_problem}) by other ways. One way will be
by multiplying some suitable constant, say $r_1\in\mathbb{N}$, on
numerator, denominator of right side of
$\frac{4}{24l+1}-\frac{1}{6l+r}=\frac{4r-1}{(6l+r)(24l+1)}$ such that
$\frac{4}{24l+1}-\frac{1}{6l+r}=\frac{4r-1}{(6l+r)(24l+1)}=\frac{(4r-1)r_1}{r_1(6l+r)(24l+1)}$
and at least one of the sums of numerator $(4r-1)r_1$ in right side
of the expression, say, $a+((4r-1)r_1-a)$ for at least one of the
possible value of $r$ such that $a,(4r-1)r_1-a$ divide
$r_1(6l+r)(24l+1)$; then $\frac{4}{24l+1}$ has the expression in the
form of the equation (\ref{Erdos-Straus_problem}). However, there is
no way to find $r$, $r_1$ in this process but trial and error method. 

A suitable method to find the expression of the equation
(\ref{Erdos-Straus_problem}) for all $n\in\mathbb{N}$ has been
discussed in Section \ref{Erdos-Straus_3rd_theorem}.

\subsubsection{Example} 
\label{Erdos-Straus_problem_common_exceptional_cases}
(i) (for $l=17$)
$\frac{4}{409}-\frac{1}{104}=\frac{7}{2^3\times 13 \times
  409}=\frac{2\times 7}{2^4\times 13 \times 409}=\frac{1+13}{2^4\times
  13 \times 409}=\frac{1}{85072}+\frac{1}{6544}$
i.e. $\frac{4}{409}=\frac{1}{104}+\frac{1}{85072}+\frac{1}{6544}$. 

(ii)
(for $l=24$) $\frac{4}{577}-\frac{1}{145}=\frac{3}{5\times 29 \times
  577}=\frac{2\times 3}{2\times 5 \times 29\times
  577}=\frac{1+5}{2\times 5 \times 29\times
  577}=\frac{1}{167330}+\frac{1}{33466}$
i.e. $\frac{4}{577}=\frac{1}{145}+\frac{1}{167330}+\frac{1}{33466}$.

(iii) (for $l=232$) $\frac{4}{5569}-\frac{1}{1394}=\frac{7}{2\times 17
  \times 41\times 5569}=\frac{6\times 7}{6\times 2\times 17 \times
  41\times 5569}=\frac{1+41}{6\times 2\times 17 \times 41\times
  5569}=\frac{1}{46579116}+\frac{1}{1136076}$
i.e. $\frac{4}{5569}=\frac{1}{1394}+\frac{1}{46579116}+\frac{1}{1136076}$.

(iv) (for $l=400$) $\frac{4}{9601}-\frac{1}{2405}=\frac{19}{5\times 13
  \times 37\times 9601}=\frac{2\times 19}{2\times 5\times 13 \times
  37\times 9601}=\frac{1+37}{2\times 5\times 13 \times 37\times
  9601}=\frac{1}{46180810}+\frac{1}{1248130}$
i.e. $\frac{4}{9601}=\frac{1}{2405}+\frac{1}{46180810}+\frac{1}{1248130}$.

(v) (for $l=997$) $\frac{4}{23929}-\frac{1}{5984}=\frac{7}{32\times 11
  \times 17\times 23929}=\frac{3\times 7}{3\times 32\times 11 \times
  17\times 23929}=\frac{4+17}{3\times 32\times 11 \times 17\times
  23929}=\frac{1}{107393352}+\frac{1}{25269024}$
i.e. $\frac{4}{23929}=\frac{1}{5984}+\frac{1}{107393352}+\frac{1}{25269024}$.

(vi) (for $l=3477$) $\frac{4}{83449}-\frac{1}{20865}=\frac{11}{3\times
  5 \times 13\times 107\times 83449}=\frac{10\times 11}{10\times
  3\times 5\times 13 \times 107\times 83449}=\frac{3+107}{10\times
  3\times 5 \times 13\times 107\times
  83449}=\frac{1}{5803877950}+\frac{1}{162725550}$
i.e. $\frac{4}{83449}=\frac{1}{20865}+\frac{1}{5803877950}+\frac{1}{162725550}$.

(vii) (for $l=4250$)
$\frac{4}{102001}-\frac{1}{25502}=\frac{7}{2\times 41 \times 311\times
  102001}=\frac{6\times 7}{6\times 2\times 41 \times 311\times
  102001}=\frac{1+41}{6\times 2\times 41 \times 311\times
  102001}=\frac{1}{15607377012}+\frac{1}{380667732}$
i.e. $\frac{4}{102001}=\frac{1}{25502}+\frac{1}{15607377012}+\frac{1}{380667732}$.

(viii) (for $l=13734$)
$\frac{4}{329617}-\frac{1}{82405}=\frac{3}{5\times 16481\times
  329617}=\frac{2\times 3}{2\times 5\times 16481\times
  329617}=\frac{1+5}{2\times 5\times 16481\times
  329617}=\frac{1}{54324177770}+\frac{1}{10864835554}$
i.e. $\frac{4}{102001}=\frac{1}{82405}+\frac{1}{54324177770}+\frac{1}{10864835554}$.

(ix) (for $l=29680$)
$\frac{4}{712321}-\frac{1}{178086}=\frac{23}{2\times 3 \times 67\times
  443\times 712321}=\frac{3\times 23}{3\times 2\times 3 \times
  67\times 443\times 712321}=\frac{2+67}{3\times 2\times 3 \times
  67\times 443\times
  712321}=\frac{1}{190281596409}+\frac{1}{5680047654}$
i.e. $\frac{4}{712321}=\frac{1}{178086}+\frac{1}{190281596409}+\frac{1}{5680047654}$.

(x) (for $l=47260$) $\frac{4}{1134241}-\frac{1}{283561}=\frac{3}{233
  \times 1217\times 1134241}=\frac{78\times 3}{78\times 233 \times
  1217\times 1134241}=\frac{1+233}{78\times 233 \times 1217\times
  1134241}=\frac{1}{25086867951678}+\frac{1}{107668961166}$
i.e. $\frac{4}{1134241}=\frac{1}{283561}+\frac{1}{25086867951678}+\frac{1}{107668961166}$.

(xi) (for $l=71842$)
$\frac{4}{1724209}-\frac{1}{431054}=\frac{7}{2\times 13\times 59\times
  281\times 1724209}=\frac{2\times 7}{2\times 2\times 13\times
  59\times 281\times 1724209}\\=\frac{1+13}{2\times 2\times 13\times
  59\times 281\times
  1724209}=\frac{1}{1486454372572}+\frac{1}{114342644044}$
i.e. $\frac{4}{1724209}=\frac{1}{431054}+\frac{1}{1486454372572}+\frac{1}{114342644044}$.

(xii) (for $l=71925$)
$\frac{4}{1726201}-\frac{1}{431566}=\frac{63}{2\times 19\times
  41\times 277\times 1726201}=\frac{5\times 63}{5\times 2\times
  19\times 41\times 277\times 1726201}\\=\frac{38+277}{5\times 2\times 19\times 41\times 277\times
  1726201}=\frac{1}{98022323785}+\frac{1}{13447105790}$
i.e. $\frac{4}{1724209}=\frac{1}{431566}+\frac{1}{98022323785}+\frac{1}{13447105790}$.

\subsection{Lemma}
\label{Erdos-Straus_1st_lemma}
If $6l+1$ or $24l+1$ has a factor $3b+2$ ($l,b\in \mathbb{N}$), then
$\frac{4}{24l+1}$ has the expression in the form of
(\ref{Erdos-Straus_problem}).

{\bf Proof}: If $6l+1$ or $24l+1$ ($l\in \mathbb{N}$) has a factor
$3b+2$, 

then $\frac{4}{24l+1}-\frac{1}{6l+1}=\frac{3}{(3b+2)\times f}$
(where $(24l+1)(6l+1)=f\times (3b+2))$ 

$=\frac{3(b+1)}{(3b+2)\times (b+1)\times
  f}=\frac{1+(3b+2)}{(3b+2)\times (b+1)\times f}=\frac{1}{(3b+2)\times
  (b+1)\times f}+\frac{1}{(b+1)\times f}$

i.e. $\frac{4}{24l+1}=\frac{1}{6l+1}+\frac{1}{(3b+2)\times (b+1)\times
  f}+\frac{1}{(b+1)\times f}$.

Example: $\frac{4}{97}-\frac{1}{25}=\frac{3}{5^2\times 97}$ ($b=1$)
$\frac{2\times 3}{2\times 5^2\times 97}=\frac{1+5}{2\times 5^2\times
  97}=\frac{1}{4850}+\frac{1}{970}$.

\subsection{Lemma}
\label{Erdos-Straus_2nd_lemma}
If $3l+1=5b$ ($l,b\in \mathbb{N}$), then
$\frac{4}{40b-7}$ has the expression in the form of
(\ref{Erdos-Straus_problem}).

{\bf Proof}: If $3l+1=5b$ ($l,b\in \mathbb{N}$), then
$\frac{4}{24l+1}-\frac{1}{6l+2}=\frac{4}{40b-7}-\frac{1}{10b}=\frac{7}{10b(40b-7)}=\frac{2+5}{10b(40b-7)}=\frac{1}{5b(40b-7)}+\frac{1}{2b(40b-7)}$
i.e. $\frac{4}{40b-7}=\frac{1}{10b}+\frac{1}{5b(40b-7)}+\frac{1}{2b(40b-7)}$.

{\bf Note:} If $24l+1$ has a factor $5$ ($l\in \mathbb{N}$), then
$\frac{4}{24l+1}$ has the expression in the form of
(\ref{Erdos-Straus_problem}).

\subsection{Lemma}
\label{Erdos-Straus_3rd_lemma}
If $3l+1=7b$ ($l,b\in \mathbb{N}$), then
$\frac{4}{56b-7}$ has the expression in the form of
(\ref{Erdos-Straus_problem}).

{\bf Proof}: If $3l+1=7b$ ($l,b\in \mathbb{N}$), then
$\frac{4}{24l+1}-\frac{1}{6l+2}=\frac{4}{56b-7}-\frac{1}{14b}=\frac{7}{14b(56b-7)}=\frac{1}{2b(56b-7)}=\frac{1}{4b(56b-7)}+\frac{1}{4b(56b-7)}$
i.e. $\frac{4}{56b-7}=\frac{1}{14b}+\frac{1}{4b(56b-7)}+\frac{1}{4b(56b-7)}$.

{\bf Note:} If $24l+1$  has a factor $7$ ($l\in \mathbb{N}$), then
$\frac{4}{24l+1}$ has the expression in the form of
(\ref{Erdos-Straus_problem}).

\subsection{Lemma}
\label{Erdos-Straus_4th_lemma}
If $3l+1=7b+5$ ($l,b\in \mathbb{N}$), then
$\frac{4}{56b+33}$ has the expression in the form of
(\ref{Erdos-Straus_problem}).

{\bf Proof}: If $3l+1=7b+5$ ($l,b\in \mathbb{N}$), then
$\frac{4}{24l+1}-\frac{1}{6l+2}=\frac{4}{56b+33}-\frac{1}{2(7b+5)}=\frac{7}{2(7b+5)(56b+33)}=\frac{7(b+1)}{2(b+1)(7b+5)(56b+33)}=\frac{2+(7b+5)}{2(b+1)(7b+5)(56b+33)}=\frac{1}{(b+1)(7b+5)(56b+33)}+\frac{1}{2(b+1)(56b+33)}$

i.e. $\frac{4}{56b+33}=\frac{1}{2(7b+5)}+\frac{1}{(b+1)(7b+5)(56b+33)}+\frac{1}{2(b+1)(56b+33)}$.

{\bf Note:} If $24l+1$  has a factor $7b+5$ ($l,b\in \mathbb{N}$), then
$\frac{4}{24l+1}$ has the expression in the form of
(\ref{Erdos-Straus_problem}).

\subsection{Lemma}
\label{Erdos-Straus_5th_lemma}
If $3l+1=7b+6$ ($l,b\in \mathbb{N}$), then
$\frac{4}{56b+41}$ has the expression in the form of
(\ref{Erdos-Straus_problem}).

{\bf Proof}: If $3l+1=7b+6$ ($l,b\in \mathbb{N}$), then
$\frac{4}{24l+1}-\frac{1}{6l+2}=\frac{4}{56b+41}-\frac{1}{2(7b+6)}=\frac{7}{2(7b+6)(56b+41)}=\frac{7(b+1)}{2(b+1)(7b+6)(56b+41)}=\frac{1+(7b+6)}{2(b+1)(7b+6)(56b+41)}=\frac{1}{2(b+1)(7b+6)(56b+41)}+\frac{1}{2(b+1)(56b+41)}$

i.e. $\frac{4}{56b+41}=\frac{1}{2(7b+6)}+\frac{1}{2(b+1)(7b+6)(56b+41)}+\frac{1}{2(b+1)(56b+41)}$.

{\bf Note:} If $24l+1$  has a factor $7b+6$ ($l,b\in \mathbb{N}$), then
$\frac{4}{24l+1}$ has the expression in the form of
(\ref{Erdos-Straus_problem}).

\subsection{Lemma}
\label{Erdos-Straus_6th_lemma}
If $3l+1=7b+3$ ($l,b\in \mathbb{N}$), then
$\frac{4}{56b+17}$ has the expression in the form of
(\ref{Erdos-Straus_problem}).

{\bf Proof}: If $3l+1=7b+3$ ($l,b\in \mathbb{N}$), then
$\frac{4}{24l+1}-\frac{1}{6l+2}=\frac{4}{56b+17}-\frac{1}{2(7b+3)}=\frac{7}{2(7b+3)(56b+17)}=\frac{7(2b+1)}{2(2b+1)(7b+3)(56b+17)}=\frac{1+(14b+6)}{(2b+1)(14b+6)(56b+41)}=\frac{1}{2(2b+1)(7b+3)(56b+17)}+\frac{1}{(2b+1)(56b+17)}$
i.e. $\frac{4}{56b+17}=\frac{1}{2(7b+3)}+\frac{1}{2(2b+1)(7b+3)(56b+17)}+\frac{1}{(2b+1)(56b+17)}$.

{\bf Note:} If $24l+1$  has a factor $7b+3$ ($l,b\in \mathbb{N}$), then
$\frac{4}{24l+1}$ has the expression in the form of
(\ref{Erdos-Straus_problem}).

\section{Theorem}
\label{Erdos-Straus_1st_theorem}
The expression in the form of (\ref{Erdos-Straus_problem}) has a
solution for every natural number $n$, except possibly for those primes of the
form $n\equiv r$ (mod 120), with $r=1, 7^2$.

{\bf Proof}: (i) If $l=5b-4$ ($l,b\in \mathbb{N}$), then
$\frac{4}{24l+1}-\frac{1}{6l+1}=\frac{4}{120b-95}-\frac{1}{30b-23}=\frac{3}{5(24b-19)(30b-23)}=\frac{2\times
  3}{2\times
  5(24b-19)(30b-23)}=\frac{1+5}{10(24b-19)(30b-23))}=\frac{1}{10(24b-19)(30b-23)}+\frac{1}{2(24b-19)(30b-23)}$
i.e. $\frac{4}{120b-95}=\frac{1}{30b-23}+\frac{1}{10(24b-19)(30b-23)}+\frac{1}{2(24b-19)(30b-23)}$.

(ii) If $l=5b-2$ ($l,b\in \mathbb{N}$), then
$\frac{4}{24l+1}-\frac{1}{6l+2}=\frac{4}{120b-47}-\frac{1}{30b-10}=\frac{7}{10(120b-47)(3b-1)}=\frac{2+5}{10(120b-47)(3b-1))}=\frac{1}{5(120b-47)(3b-1)}+\frac{1}{2(120b-47)(3b-1)}$
i.e. $\frac{4}{120b-47}=\frac{1}{30b-10}+\frac{1}{5(120b-47)(3b-1)}+\frac{1}{2(120b-47)(3b-1)}$.

(iii) If $l=5b-1$ ($l,b\in \mathbb{N}$), then
$\frac{4}{24l+1}-\frac{1}{6l+1}=\frac{4}{120b-23}-\frac{1}{30b-5}=\frac{3}{5(120b-23)(6b-1)}=\frac{2\times
  3}{2\times
  5(120b-23)(6b-1)}=\frac{1+5}{10(120b-23)(6b-1)}=\frac{1}{10(120b-23)(6b-1)}+\frac{1}{2(120b-23)(6b-1)}$
i.e. $\frac{4}{120b-23}=\frac{1}{30b-5}+\frac{1}{10(120b-23)(6b-1)}+\frac{1}{2(120b-23)(6b-1)}$.

(iv) If $l=5b-3$ ($l,b\in \mathbb{N}$), then
$\frac{4}{24l+1}=\frac{4}{120b-71}$ where $120b-71=49$ (mod 120) 

(v) If $l=5b$ ($l,b\in \mathbb{N}$), then $\frac{4}{24l+1}=\frac{4}{120b+1}$
where $120b+1=1$ (mod 120). 

Hence, the expression in the form of (\ref{Erdos-Straus_problem}) has
a solution for every natural number $n$, except possibly for those
primes of the form $n\equiv r$ (mod 120), with $r=1, 7^2$.

\section{Theorem (Mordell)}
\label{Erdos-Straus_2nd_theorem}
The expression in the form of (\ref{Erdos-Straus_problem}) has a
solution for every number $n$, except possibly for those primes of the
form $n\equiv r$ (mod 780), with $r=1^2,11^2,13^2,17^2,19^2,23^2$.

{\bf Proof}: From Theorem \ref{Erdos-Straus_1st_theorem}, we know that
the expression (\ref{Erdos-Straus_problem}) has a solution if
$\frac{4}{n}=\frac{4}{120b+r}$ except $r=1,-71$.

(1) Let $\frac{4}{n}=\frac{4}{120b+1}$.

(i) If $b=7c-6$ ($b,c\in \mathbb{N}$), then $\frac{4}{120b+1}=\frac{4}{840c-719}$ where $840c-719=11^2$ (mod 840).

(ii) If $b=7c-5$ ($b,c\in \mathbb{N}$), then
$\frac{4}{120b+1}-\frac{1}{30b+3}=\frac{4}{840c-599}-\frac{1}{210c-147}=\frac{11}{21(10c-7)(840c-599)}=\frac{2\times
  11}{2\times
  21(10c-7)(840c-599)}=\frac{1+21}{42(10c-7)(840c-599)}=\frac{1}{42(10c-7)(840c-599)}+\frac{1}{2(10c-7)(840c-599)}$
i.e. $\frac{4}{840c-599}=\frac{1}{210c-147}+\frac{1}{42(10c-7)(840c-599)}+\frac{1}{2(10c-7)(840c-599)}$.

(iii) If $b=7c-4$ ($b,c\in \mathbb{N}$), then
$\frac{4}{120b+1}=\frac{4}{840c-479}$ where $840c-479=19^2$ (mod 840).

(iv) If $b=7c-3$ ($b,c\in \mathbb{N}$), then
$\frac{4}{120b+1}-\frac{1}{30b+2}=\frac{4}{840c-359}-\frac{1}{210c-88}=\frac{7}{2(105c-44)(840c-359)}=\frac{(15c-6)\times
  7}{(15c-6)\times
  2(105c-44)(840c-359)}=\frac{2+(105c-42)}{2(15c-6)(105c-44)(840c-359)}=\frac{1}{(15c-6)(105c-44)(840c-359)}+\frac{1}{2(15c-6)(840c-359)}$
i.e. $\frac{4}{840c-359}=\frac{1}{210c-88}+\frac{1}{(15c-6)(105c-44)(840c-359)}+\frac{1}{2(15c-6)(840c-359)}$.

(v) If $b=7c-2$ ($b,c\in \mathbb{N}$), then
$\frac{4}{120b+1}-\frac{1}{30b+2}=\frac{4}{840c-239}-\frac{1}{210c-58}=\frac{7}{2(105c-29)(840c-239)}=\frac{(15c-4)\times
  7}{(15c-4)\times
  2(105c-29)(840c-239)}=\frac{1+(105c-29)}{2(15c-4)(105c-29)(840c-239)}=\frac{1}{2(15c-4)(105c-29)(840c-239)}+\frac{1}{2(15c-4)(840c-239)}$
i.e. $\frac{4}{840c-239}=\frac{1}{210c-58}+\frac{1}{2(15c-4)(105c-29)(840c-239)}+\frac{1}{2(15c-4)(840c-239)}$.

(vi) If $b=7c-1$ ($b,c\in \mathbb{N}$), then
$\frac{4}{120b+1}-\frac{1}{30b+2}=\frac{4}{840c-119}-\frac{1}{210c-28}=\frac{7}{7(120c-17)14(15c-2)}=\frac{1}{14(15c-2)(120c-17)}=\frac{1}{28(15c-2)(120c-17)}+\frac{1}{28(15c-2)(120c-17)}$
i.e. $\frac{4}{840c-119}=\frac{1}{210c-28}+\frac{1}{28(15c-2)(120c-17)}+\frac{1}{28(15c-2)(120c-17)}$.

(vii) If $b=7c$ ($b,c\in \mathbb{N}$), then
$\frac{4}{120b+1}=\frac{4}{840c+1}$ where $840c+1=1^2$ (mod 840).

(2) Let $\frac{4}{n}=\frac{4}{120b-71}$.

(i) If $b=7c-6$ ($b,c\in \mathbb{N}$), then
$\frac{4}{120b-71}-\frac{1}{30b-16}=\frac{4}{840c-791}-\frac{1}{210c-196}=\frac{7}{7(120c-113)14(15c-14)}=\frac{1}{14(120c-113)(15c-14)}=\frac{1}{28(120c-113)(15c-14)}+\frac{1}{28(120c-113)(15c-14)}$
i.e. $\frac{4}{840c-791}=\frac{1}{210c-196}+\frac{1}{28(120c-113)(15c-14)}+\frac{1}{28(120c-113)(15c-14)}$.

(ii) If $b=7c-5$ ($b,c\in \mathbb{N}$), then $\frac{4}{120b-71}=\frac{4}{840c-671}$ where $840c-671=13^2$ (mod 840).

(iii) If $b=7c-4$ ($b,c\in \mathbb{N}$), then $\frac{4}{120b-71}=\frac{4}{840c-551}$ where $840c-551=17^2$ (mod 840).

(iv) If $b=7c-3$ ($b,c\in \mathbb{N}$), then
$\frac{4}{120b-71}-\frac{1}{30b-16}=\frac{4}{840c-431}-\frac{1}{210c-106}=\frac{7}{2(105c-53)(840c-431)}=\frac{(30c-15)\times
  7}{(30c-15)\times
  2(105c-53)(840c-431)}=\frac{1+(210c-106)}{(30c-15)(210c-106)(840c-431)}=\frac{1}{30(2c-1)(105c-53)(840c-431)}+\frac{1}{15(2c-1)(840c-431)}$
i.e. $\frac{4}{840c-431}=\frac{1}{210c-106}+\frac{1}{30(2c-1)(105c-53)(840c-431)}+\frac{1}{15(2c-1)(840c-431)}$.

(v) If $b=7c-2$ ($b,c\in \mathbb{N}$), then $\frac{4}{120b-71}=\frac{4}{840c-311}$ where $840c-311=23^2$ (mod 840).

(vi) If $b=7c-1$ ($b,c\in \mathbb{N}$), then
$\frac{4}{120b-71}-\frac{1}{30b-16}=\frac{4}{840c-191}-\frac{1}{210c-46}=\frac{7}{2(105c-23)(840c-191)}=\frac{(15c-3)\times
  7}{(15c-3)\times
  2(105c-23)(840c-191)}=\frac{2+(105c-23)}{2(15c-3)(105c-23)(840c-191)}=\frac{1}{3(5c-1)(105c-23)(840c-191)}+\frac{1}{6(5c-1)(840c-191)}$
i.e. $\frac{4}{840c-191}=\frac{1}{210c-46}+\frac{1}{3(5c-1)(105c-23)(840c-191)}+\frac{1}{6(5c-1)(840c-191)}$.

(vii) If $b=7c$ ($b,c\in \mathbb{N}$), then
$\frac{4}{120b-71}-\frac{1}{30b-16}=\frac{4}{840c-71}-\frac{1}{210c-16}=\frac{7}{2(105c-8)(840c-71)}=\frac{(15c-1)\times
  7}{(15c-1)\times
  2(105c-8)(840c-71)}=\frac{1+(105c-8)}{2(15c-1)(105c-8)(840c-71)}=\frac{1}{2(15c-1)(105c-8)(840c-71)}+\frac{1}{2(15c-1)(840c-71)}$
i.e. $\frac{4}{840c-71}=\frac{1}{210c-16}+\frac{1}{2(15c-1)(105c-8)(840c-71)}+\frac{1}{2(15c-1)(840c-71)}$.

Hence, the expression in the form of the equation
(\ref{Erdos-Straus_problem}) has a solution for every number $n$,
except possibly for those primes of the form $n\equiv r$ (mod 840),
with $r=1^2,11^2,13^2,17^2,19^2,23^2$.

\section{Theorem}
\label{Erdos-Straus_3rd_theorem}
The Erd\H{o}s-Straus conjecture is true i.e. the expression in the
form of (\ref{Erdos-Straus_problem}) has a solution for every natural
number $n$.

{\bf Proof}: We know that the Erd\H{o}s-Straus conjecture is
equivalent to find the equation (\ref{Erdos-Straus_problem}) for all
primes. Because, it is well known that if
$\frac{4}{n}=\frac{1}{x}+\frac{1}{y}+\frac{1}{z}$ is true for $n=p$,
then we get $\frac{4}{mp}=\frac{1}{mx}+\frac{1}{my}+\frac{1}{mz}$ for
any natural number $m$. Let us consider, $n$ as prime number $p$,
which is greater than 2, in the equation
(\ref{Erdos-Straus_problem}). For $p=2$, we know
$\frac{4}{2}=\frac{1}{1}+\frac{1}{2}+\frac{1}{2}$.

{\bf Case 1}: If $x=y=z$, then $3p=4x$ which is
a invalid equation. Thus the equation
(\ref{Erdos-Straus_problem}) has no solution if $n$ is a prime number and $x=y=z$.

{\bf Case 2}: Two of $x,y,z$ are equal. Without loss of generality,
let $x=y$. 
\begin{equation}
\text{Then}~\frac{4}{p}=\frac{2}{x}+\frac{1}{z}~\text{i.e.}~p(x+2z)=4xz. 
\label{Erdos-Straus_problem_x_and_y_equal}
\end{equation}
Thus $p$ divides $x$ or $p$ divides $y$ for $p>2$.

(I) Let $p$ divides $x$, then $x=up,~u\in \mathbb{N}$. From
(\ref{Erdos-Straus_problem_x_and_y_equal}), we get $up+2z=4uz$. Or, $2z(2u-1)=up$.

(i) If $p$ divides $z$ i.e. $z=vp,~v\in \mathbb{N}$, then
$2v(2u-1)=u$, which is a invalid equation.

(ii) If $p$ divides $2u-1$ i.e. $2u-1=v_1p,~v_1\in \mathbb{N}$, then
$2zv_1=u$ and $u$ is even since $z=\frac{u}{2v_1}$ i.e. $2v_1$ divides
$u$. Then $u=2u_1$, $4u_1-1=v_1p$ and $v_1$ divides $u_1$. So, $v_1$ divides
$u_1$ and $v_1$ divides $4u_1-1$ imply that $v_1=1$. Thus
$p=4u_1-1=4u_2+3,~x=2(u_2+1)p=y,~z=u_2+1$, where $u_1=u_2+1$.
\begin{equation}
\text{Then}~\frac{4}{p}=\frac{1}{2(u_2+1)p}+\frac{1}{2(u_2+1)p}+\frac{1}{u_2+1}=\frac{1}{\frac{p(p+1)}{2}}+\frac{1}{\frac{p(p+1)}{2}}+\frac{1}{\frac{p+1}{4}}. 
\label{Erdos-Straus_problem_p_4n_and_3_1st}
\end{equation}

(II) Let $p$ divides $z$, then $z=u_3p,~u_3\in \mathbb{N}$. From
(\ref{Erdos-Straus_problem_x_and_y_equal}), we get $x+2u_3p=4u_3x$. Or, $(4u_3-1)x=2u_3p$.

(i) If $p$ divides $x$ i.e. $x=v_2p,~v_2\in \mathbb{N}$, then
$(4u_3-1)v_2=2u_3$, which is a invalid equation.

(ii) If $p$ divides $4u_3-1$ i.e. $4u_3-1=v_3p,~v_3\in \mathbb{N}$, then
$xv_3=2u_3$. Thus $v_3$ divides
$2u_3$ and $v_3$ divides $4u_3-1$ imply that $v_3=1$. Thus
$p=4u_3-1=4u_4+3,~x=2(u_4+1)=y,~z=(u_4+1)p$, where $u_3=u_4+1$.
\begin{equation}
\text{Then}~\frac{4}{p}=\frac{1}{2(u_4+1)}+\frac{1}{2(u_4+1)}+\frac{1}{(u_2+1)p}=\frac{1}{\frac{(p+1)}{2}}+\frac{1}{\frac{(p+1)}{2}}+\frac{1}{\frac{p(p+1)}{4}}. 
\label{Erdos-Straus_problem_p_4n_and_3_2nd}
\end{equation}

{\bf Case 3}: Let $x\neq y\neq z \neq x$. From the equation
(\ref{Erdos-Straus_problem}), we get $p(xy+yz+zx)=4xyz$. Then $p$
divides at least one of $x,y,z$. Without loss of generality, let $p$
divides $x$ i.e. $x=u_5p,~u_5\in \mathbb{N}$. 
\begin{equation}
\text{Thus}~u_5p(y+z)=yz(4u_5-1). 
\label{Erdos-Straus_problem_case3_1st_eqn}
\end{equation}

(I) Let $p$ divides $y$ or $z$. Without loss of generality, let $p$
divides $y$. Then $y=v_4p,~v_4\in \mathbb{N}$. From
(\ref{Erdos-Straus_problem_case3_1st_eqn}), we get $u_5v_4p=z(4u_5v_4-u_5-v_4)$. 

(i) If $p$ divides $z$ i.e. $z=w_1p,~w_1\in \mathbb{N}$, then
$u_5v_4+w_1v_4+w_1u_5=4u_5v_4w_1$, which is a invalid equation.

(ii) If $p$ divides $4u_5v_4-u_5-v_4$,
\begin{equation}
\text{then}~4u_5v_4-u_5-v_4=w_2p,~w_2\in \mathbb{N}. 
\label{Erdos-Straus_problem_case3_2nd_eqn}
\end{equation}
So, $u_5v_4=zw_2$. Thus $w_2$ divides $u_5v_4$, $w_2$ divides
$4u_5v_4-u_5-v_4$ and hence $w_2$ divide $u_5+v_4$, $u_5^2$, $v_4^2$,
$u_5(u_5-v_4)$, $v_4(u_5-v_4)$, $(u_5-v_4)^2$ etc. Now, $w_2$ divide
$u_5v_4$ and $u_5+v_4$ mean $u_5v_4=w_3w_2,~w_3\in \mathbb{N}$ and
$u_5+v_4=w_4w_2,~w_4\in \mathbb{N}$. Then from
(\ref{Erdos-Straus_problem_case3_2nd_eqn}), we get $p=4w_3-w_4$. Hence
the equation (\ref{Erdos-Straus_problem}) has solution when
$n=p=4w_3-w_4$ with $x=u_5p$, $y=v_4p$, $z=w_3$, $u_5v_4=w_3w_2$,
$u_5+v_4=w_4w_2$ since
\begin{equation}
\frac{1}{x}+\frac{1}{y}+\frac{1}{z}=\frac{(u_5+v_4)w_3+u_5v_4p}{u_5v_4w_3p}=\frac{w_4w_2w_3+w_3w_2p}{u_5v_4w_3p}=\frac{(w_4+p)w_2}{u_5v_4p}=\frac{4w_3w_2}{w_3w_2p}=\frac{4}{p}. 
\label{Erdos-Straus_problem_case3_3rd_eqn}
\end{equation}
We have to demonstrate the solutions of the equation
(\ref{Erdos-Straus_problem}) if $n=4m+1$ since we already know
solutions of it for $n=4m$, $n=4m+2$ and $n=4m+3$.

(a) If $w_4=3$, then $p=4w_3-3$ (i.e. $n=4m+1$), $v_4=3w_2-u_5$,
$w_3w_2=u_5v_4=u_5(3w_2-u_5)$
i.e. $w_2=\frac{u_5^2}{3u_5-w_3}=\frac{u_5^2}{3u_5-\frac{p+3}{4}}$. Thus
we have to choose $u_5$
$\left(>\left[\frac{\left(\frac{p+3}{4}\right)}{3}\right]\right)$ such
that $w_2\in \mathbb{N}$. However, it will have limited use for the
solutions when $n=4m+1$ as $w_3=3$ is a special case.

For example, if $p=13$, then $w_3=4,~w_4=3$. If we choose $u_5=2$,
then $w_2=2,~v_4=4$. Hence $\frac{4}{13}=\frac{1}{2\times
  13}+\frac{1}{4\times 13}+\frac{1}{4}$. We can try to reproduce the solutions
of Section \ref{Erdos-Straus_problem_common_exceptional_cases}.

For $p=409$, $w_3=103,~w_4=3$. Then, we can not find any $w_2\in
\mathbb{N}$ if $u_5<1000001$.

For $p=577$, $w_3=145,~w_4=3$. If we choose $u_5=50$, then we can get
$w_2=500,~v_4=1450$. Thus $\frac{4}{577}=\frac{1}{50\times
  577}+\frac{1}{1450\times 577}+\frac{1}{145}$. If we choose $u_5=58$,
then we can get $w_2=116,~v_4=290$ and
$\frac{4}{577}=\frac{1}{58\times 577}+\frac{1}{290\times
  577}+\frac{1}{145}$.

For $p=5569$, $w_3=1393,~w_4=3$. Then, we can not find any $w_2\in
\mathbb{N}$ if $u_5<1000001$.

For $p=9601$, $w_3=2401,~w_4=3$. Then, we can not find any $w_2\in
\mathbb{N}$ if $u_5<1000001$.

{\bf Note}: Let $w_2\in \mathbb{N}$, then $3u_5-w_3>0$. Thus
$v_4=3w_2-u_5=\frac{3u_5^2}{3u_5-w_3}-u_5=\frac{u_5w_3}{3u_5-w_3}>0$. Hence
$w_2\in \mathbb{N}$ implies $v_4\in \mathbb{N}$.

(b) Let $w_4=4w_5+3$ where $w_5=0$ or $w_5\in \mathbb{N}$, then
$p=4(w_3-w_5)-3$ (i.e. $n=p=4m+1$). So, $w_3-w_5=\frac{p+3}{4}\in
\mathbb{N}$. Or, $w_3=w_5+\frac{p+3}{4}$,
$w_2=\frac{u_5^2}{(4w_5+3)u_5-\frac{p+3}{4}-w_5}$. Thus, we have to
choose $w_5$, $u_5$
$\left(>\left[\frac{w_5+\frac{p+3}{4}}{4w_5+3}\right]\right)$ such
that $w_2\in \mathbb{N}$. Using these $w_5$, $u_5$, $w_2$; we can
calculate $w_3=w_5+\frac{p+3}{4}$, $v_4=(4w_5+3)w_2-u_5$. Hence, in
this way we can generate all solutions of the equation
(\ref{Erdos-Straus_problem}) if $n=p=4m+1$.

As examples, we can reproduce the solutions of Section
\ref{Erdos-Straus_problem_common_exceptional_cases}. We can get
multiple solutions if we choose $w_5,~u_5\le 1000$ with
$u_5>\left[\frac{w_5+\frac{p+3}{4}}{4w_5+3}\right]$, however for same
solutions as of Section
\ref{Erdos-Straus_problem_common_exceptional_cases}, we have to
increase the range of $u_5$ for majority of cases.

Let $p=409$. For $w_5,~u_5\le 1000$, we get eleven solutions of the
equation (\ref{Erdos-Straus_problem}). The first solution is for
$w_5=1,~u_5=15$. Then $w_2=225$, $v_4=1560$, $w_3=104$ and hence
$\frac{4}{409}=\frac{1}{15\times 409}+\frac{1}{1560\times
  409}+\frac{1}{104}$. The last solution is for
$w_5=14,~u_5=234$. Then $w_2=4$, $v_4=2$, $w_3=117$ and hence
$\frac{4}{409}=\frac{1}{234\times 409}+\frac{1}{2\times
  409}+\frac{1}{117}$. The same solutions as of Section
\ref{Erdos-Straus_problem_common_exceptional_cases} is for
$w_5=1,~u_5=16$. Then $w_2=32$, $v_4=208$, $w_3=104$ and hence
$\frac{4}{409}=\frac{1}{16\times 409}+\frac{1}{208\times
  409}+\frac{1}{104}$.

Let $p=577$. For $w_5,~u_5\le 1000$, we get twelve solutions of the
equation (\ref{Erdos-Straus_problem}). The first solution is for
$w_5=0,~u_5=50$. Then $w_2=500$, $v_4=1450$, $w_3=145$ and hence
$\frac{4}{577}=\frac{1}{50\times 577}+\frac{1}{1450\times
  577}+\frac{1}{145}$. The last solution is for
$w_5=20,~u_5=330$. Then $w_2=4$, $v_4=2$, $w_3=165$ and hence
$\frac{4}{577}=\frac{1}{330\times 577}+\frac{1}{2\times
  577}+\frac{1}{165}$. The same solutions as of Section
\ref{Erdos-Straus_problem_common_exceptional_cases} is for
$w_5=0,~u_5=58$. Then $w_2=116$, $v_4=290$, $w_3=104$ and hence
$\frac{4}{577}=\frac{1}{58\times 577}+\frac{1}{290\times
  577}+\frac{1}{104}$.

Let $p=5569$. For $w_5,~u_5\le 1000$, we get eleven solutions of the
equation (\ref{Erdos-Straus_problem}). The first solution is for
$w_5=1,~u_5=204$. Then $w_2=1224$, $v_4=8364$, $w_3=1394$ and hence
$\frac{4}{5569}=\frac{1}{204\times 5569}+\frac{1}{8364\times
  5569}+\frac{1}{1394}$, which is the same solutions as of Section
\ref{Erdos-Straus_problem_common_exceptional_cases}. The last solution
is for $w_5=35,~u_5=10$. Then $w_2=50$, $v_4=7140$, $w_3=1428$ and
hence $\frac{4}{5569}=\frac{1}{10\times 5569}+\frac{1}{7140\times
  5569}+\frac{1}{1428}$. 

Let $p=9601$. For $w_5,~u_5\le 1000$, we get six solutions of the
equation (\ref{Erdos-Straus_problem}). The first solution is for
$w_5=4,~u_5=130$. Then $w_2=260$, $v_4=4810$, $w_3=2405$ and hence
$\frac{4}{9601}=\frac{1}{130\times 9601}+\frac{1}{4810\times
  9601}+\frac{1}{2405}$, which is the same solutions as of Section
\ref{Erdos-Straus_problem_common_exceptional_cases}. The last solution
is for $w_5=104,~u_5=6$. Then $w_2=4$, $v_4=1670$, $w_3=2505$ and
hence $\frac{4}{9601}=\frac{1}{6\times 9601}+\frac{1}{1670\times
  9601}+\frac{1}{2505}$.

Let $p=23929$. For $w_5,~u_5\le 1000$, we get twenty four solutions of
the equation (\ref{Erdos-Straus_problem}). The first solution is for
$w_5=1,~u_5=855$. Then $w_2=731025$, $v_4=5116320$, $w_3=5984$ and
hence $\frac{4}{23929}=\frac{1}{855\times
  23929}+\frac{1}{5116320\times 23929}+\frac{1}{5984}$. The last
solution is for $w_5=854,~u_5=2$. Then $w_2=4$, $v_4=13674$,
$w_3=6837$ and hence $\frac{4}{23929}=\frac{1}{2\times
  23929}+\frac{1}{13674\times 23929}+\frac{1}{6837}$. To get the same
solutions as of Section
\ref{Erdos-Straus_problem_common_exceptional_cases}, we can consider
the case $w_5\le 1000,~u_5\le 10000$ and then we get thirty six
solutions of the equation (\ref{Erdos-Straus_problem}). Thus the same
solutions as of Section
\ref{Erdos-Straus_problem_common_exceptional_cases} is for
$w_5=1,~u_5=1056$. Then $w_2=792$, $v_4=4488$, $w_3=5984$ and hence
$\frac{4}{23929}=\frac{1}{1056\times 23929}+\frac{1}{4488\times
  23929}+\frac{1}{5984}$.

Let $p=83449$. For $w_5,~u_5\le 1000$, we get eleven solutions of
the equation (\ref{Erdos-Straus_problem}). The first solution is for
$w_5=5,~u_5=908$. Then $w_2=51529$, $v_4=1184259$, $w_3=20868$ and
hence $\frac{4}{83449}=\frac{1}{908\times
  83449}+\frac{1}{1184259\times 83449}+\frac{1}{20868}$. The last
solution is for $w_5=353,~u_5=15$. Then $w_2=25$, $v_4=35360$,
$w_3=21216$ and hence $\frac{4}{83449}=\frac{1}{15\times
  83449}+\frac{1}{35360\times 83449}+\frac{1}{21216}$. To get the same
solutions as of Section
\ref{Erdos-Straus_problem_common_exceptional_cases}, we can consider
the case $w_5\le 1000,~u_5\le 10000$ and then we get seventeen
solutions of the equation (\ref{Erdos-Straus_problem}). Thus the same
solutions as of Section
\ref{Erdos-Straus_problem_common_exceptional_cases} is for
$w_5=2,~u_5=1950$. Then $w_2=6500$, $v_4=69550$, $w_3=20865$ and hence
$\frac{4}{83449}=\frac{1}{1950\times 83449}+\frac{1}{69550\times
  83449}+\frac{1}{20865}$.

Let $p=102001$. For $w_5,~u_5\le 1000$, we get thirteen solutions of
the equation (\ref{Erdos-Straus_problem}). The first solution is for
$w_5=7,~u_5=826$. Then $w_2=6962$, $v_4=214996$, $w_3=25508$ and hence
$\frac{4}{102001}=\frac{1}{826\times 102001}+\frac{1}{214996\times
  102001}+\frac{1}{25508}$. The last solution is for
$w_5=542,~u_5=12$. Then $w_2=16$, $v_4=34724$, $w_3=26043$ and hence
$\frac{4}{102001}=\frac{1}{12\times 102001}+\frac{1}{34724\times
  102001}+\frac{1}{26043}$. To get the same solutions as of Section
\ref{Erdos-Straus_problem_common_exceptional_cases}, we can consider
the case $w_5\le 1000,~u_5\le 10000$ and then we get twenty solutions
of the equation (\ref{Erdos-Straus_problem}). Thus the same solutions
as of Section \ref{Erdos-Straus_problem_common_exceptional_cases} is
for $w_5=1,~u_5=3732$. Then $w_2=22392$, $v_4=153012$, $w_3=25502$ and hence
$\frac{4}{102001}=\frac{1}{3732\times 102001}+\frac{1}{153012\times
  102001}+\frac{1}{25502}$.

Let $p=329617$. For $w_5,~u_5\le 1000$, we get fourteen solutions of
the equation (\ref{Erdos-Straus_problem}). The first solution is for
$w_5=26,~u_5=774$. Then $w_2=1548$, $v_4=164862$, $w_3=82431$ and hence
$\frac{4}{329617}=\frac{1}{774\times 329617}+\frac{1}{164862\times
  329617}+\frac{1}{82431}$. The last solution is for $w_5=867,~u_5=24$. Then
$w_2=18$, $v_4=62454$, $w_3=83272$ and hence $\frac{4}{329617}=\frac{1}{24\times
  329617}+\frac{1}{62454\times 329617}+\frac{1}{83272}$. To get the same
solutions as of Section
\ref{Erdos-Straus_problem_common_exceptional_cases}, we can consider
the case $w_5\le 1000,~u_5\le 100000$ and then we get thirty two
solutions of the equation (\ref{Erdos-Straus_problem}). Thus the same
solutions as of Section
\ref{Erdos-Straus_problem_common_exceptional_cases} is for
$w_5=0,~u_5=32962$. Then $w_2=65924$, $v_4=164810$, $w_3=82405$ and hence
$\frac{4}{329617}=\frac{1}{32962\times 329617}+\frac{1}{164810\times
  329617}+\frac{1}{82405}$.

Let $p=712321$. For $w_5,~u_5\le 1000$, we get fourteen solutions of
the equation (\ref{Erdos-Straus_problem}). The first solution is for
$w_5=47,~u_5=936$. Then $w_2=1352$, $v_4=257296$, $w_3=178128$ and
hence $\frac{4}{712321}=\frac{1}{936\times
  712321}+\frac{1}{257296\times 712321}+\frac{1}{178128}$. The last
solution is for $w_5=587,~u_5=76$. Then $w_2=722$, $v_4=1697346$,
$w_3=178668$ and hence $\frac{4}{712321}=\frac{1}{76\times
  712321}+\frac{1}{1697346\times 712321}+\frac{1}{178668}$. To get the same
solutions as of Section
\ref{Erdos-Straus_problem_common_exceptional_cases}, we can consider
the case $w_5\le 1000,~u_5\le 10000$ and then we get twenty six
solutions of the equation (\ref{Erdos-Straus_problem}). Thus the same
solutions as of Section
\ref{Erdos-Straus_problem_common_exceptional_cases} is for
$w_5=5,~u_5=7974$. Then $w_2=11961$, $v_4=267129$, $w_3=178086$ and hence
$\frac{4}{712321}=\frac{1}{7974\times 712321}+\frac{1}{267129\times
  712321}+\frac{1}{178086}$.

Let $p=1134241$. For $w_5,~u_5\le 1000$, we get twelve solutions of
the equation (\ref{Erdos-Straus_problem}). The first solution is for
$w_5=101,~u_5=697$. Then $w_2=28577$, $v_4=11630142$, $w_3=283662$ and
hence $\frac{4}{1134241}=\frac{1}{697\times
  1134241}+\frac{1}{11630142\times 1134241}+\frac{1}{283662}$. The
last solution is for $w_5=696,~u_5=102$. Then $w_2=612$,
$v_4=1705542$, $w_3=284257$ and hence
$\frac{4}{1134241}=\frac{1}{102\times 1134241}+\frac{1}{1705542\times
  1134241}+\frac{1}{284257}$. To get the same solutions as of Section
\ref{Erdos-Straus_problem_common_exceptional_cases}, we can consider
the case $w_5\le 1000,~u_5\le 100000$ and then we get forty four
solutions of the equation (\ref{Erdos-Straus_problem}). Thus the same
solutions as of Section
\ref{Erdos-Straus_problem_common_exceptional_cases} is for
$w_5=0,~u_5=94926$. Then $w_2=7404228$, $v_4=22117758$, $w_3=283561$
and hence $\frac{4}{1134241}=\frac{1}{94926\times
  1134241}+\frac{1}{22117758\times 1134241}+\frac{1}{283561}$.

Let $p=1724209$. For $w_5,~u_5\le 1000$, we get seven solutions of the
equation (\ref{Erdos-Straus_problem}). The first solution is for
$w_5=125,~u_5=858$. Then $w_2=1859$, $v_4=934219$, $w_3=431178$ and
hence $\frac{4}{1724209}=\frac{1}{858\times
  1724209}+\frac{1}{934219\times 1724209}+\frac{1}{431178}$. The last
solution is for $w_5=749,~u_5=144$. Then $w_2=384$, $v_4=1151472$,
$w_3=431802$ and hence $\frac{4}{1724209}=\frac{1}{144\times
  1724209}+\frac{1}{1151472\times 1724209}+\frac{1}{431802}$. To get
the same solutions as of Section
\ref{Erdos-Straus_problem_common_exceptional_cases}, we can consider
the case $w_5\le 1000,~u_5\le 100000$ and then we get thirty five
solutions of the equation (\ref{Erdos-Straus_problem}). Thus the same
solutions as of Section
\ref{Erdos-Straus_problem_common_exceptional_cases} is for
$w_5=1,~u_5=66316$. Then $w_2=132632$, $v_4=862108$, $w_3=431054$ and
hence $\frac{4}{1724209}=\frac{1}{66316\times
  1724209}+\frac{1}{862108\times 1724209}+\frac{1}{431054}$.

Let $p=1726201$. For $w_5,~u_5\le 1000$, we get thirteen solutions of
the equation (\ref{Erdos-Straus_problem}). The first solution is for
$w_5=125,~u_5=864$. Then $w_2=256$, $v_4=127904$, $w_3=431676$ and
hence $\frac{4}{1726201}=\frac{1}{864\times
  1726201}+\frac{1}{127904\times 1726201}+\frac{1}{431676}$. The last
solution is for $w_5=473,~u_5=228$. Then $w_2=1444$, $v_4=2736152$,
$w_3=432024$ and hence $\frac{4}{1726201}=\frac{1}{228\times
  1726201}+\frac{1}{2736152\times 1726201}+\frac{1}{432024}$. To get
the same solutions as of Section
\ref{Erdos-Straus_problem_common_exceptional_cases}, we can consider
the case $w_5\le 1000,~u_5\le 10000$ and then we get twenty seven
solutions of the equation (\ref{Erdos-Straus_problem}). Thus the same
solutions as of Section
\ref{Erdos-Straus_problem_common_exceptional_cases} is for
$w_5=15,~u_5=7790$. Then $w_2=1025$, $v_4=56785$, $w_3=431566$ and
hence $\frac{4}{1726201}=\frac{1}{7790\times
  1726201}+\frac{1}{56785\times 1726201}+\frac{1}{431566}$.

Thus the Erd\H{o}s-Straus conjecture is true for all values of natural
numbers $n\ge 2$.

{\bf Note}: Let $w_2\in \mathbb{N}$, then
$(4w_5+3)u_5-\frac{p+3}{4}-w_5>0$. Thus
$v_4=(4w_5+3)w_2-u_5=\frac{(4w_5+3)u_5^2}{(4w_5+3)u_5-\frac{p+3}{4}-w_5}-u_5=\frac{u_5(\frac{p+3}{4}+w_5)}{(4w_5+3)u_5-\frac{p+3}{4}-w_5}>0$. Hence
$w_2\in \mathbb{N}$ implies $v_4\in \mathbb{N}$.

\subsection{Remark}
For $w_4=4w_5+3$, $p=4(w_3-w_5)-3$,
$x=u_5p,~y=v_4p,~z=w_3,~u_5v_4=w_2w_3,~u_5+v_4=w_2w_4$. Then
$w_3=\frac{p+3}{4}+w_5$,
$u_5v_4=w_2(\frac{p+3}{4}+w_5),~u_5+v_4=w_2(4w_5+3)$. Thus
$\frac{u_5v_4}{u_5+v_4}=\frac{\frac{p+3}{4}+w_5}{4w_5+3}$. Hence
$u_5v_4(4w_5+3)=(\frac{p+3}{4}+w_5)(u_5+v_4)$.

\subsection{Corollary}
\label{Erdos-Straus_1st_corollary}
If $w_2=1$ and $u_5+v_4=4w_6-1$, then
\begin{equation}
p=4\left\{(u_5(4w_6-u_5-1)-w_6)\right\}+1~\text{and}~\frac{4}{p}=\frac{1}{u_5p}+\frac{1}{(4w_6-u_5-1)p}+\frac{1}{u_5(4w_6-u_5-1)}. 
\label{Erdos-Straus_problem_1st_corollary_1eq}
\end{equation}

{\bf Proof}: Let $w_2=1$. Then from equation
(\ref{Erdos-Straus_problem_case3_2nd_eqn}), we get
$p=4u_5v_4-u_5-v_4,~x=u_5p,~y=v_4p,~z=u_5v_4$. If $u_5+v_4=4w_6-1$,
then
$p=4\left\{(u_5(4w_6-u_5-1)-w_6)\right\}+1,~y=(4w_6-u_5-1)p$. Thus
$\frac{1}{u_5p}+\frac{1}{(4w_6-u_5-1)p}+\frac{1}{u_5(4w_6-u_5-1)}=\frac{(4w_6-u_5-1)+u_5+p}{u_5(4w_6-u_5-1)p}=\frac{(4w_6-1)+4\left\{(u_5(4w_6-u_5-1)-w_6)\right\}+1}{u_5(4w_6-u_5-1)p}=\frac{4u_5(4w_6-u_5-1)}{u_5(4w_6-u_5-1)p}=\frac{4}{p}$.

If (a) $u_5=1$, then $p=12w_6-7=12w_7+5$ (b) $u_5=2$, then
$p=28w_6-23=28w_7+5$ (c) $u_5=3$, then $p=44w_6-47=12w_7+41$ (d)
$u_5=4$, then $p=60w_6-79=60w_7+41$ (e) $u_5=5$, then
$p=76w_6-119=76w_7+33$ (f) $u_5=6$, then $p=92w_6-167=12w_7+17$ (g)
$u_5=7$, then $p=108w_6-223=108w_7+101$ (h) $u_5=8$, then
$p=124w_6-287=124w_7+85$ (i) $u_5=9$, then $p=140w_6-359=140w_7+61$
(j) $u_5=10$, then $p=156w_6-439=156w_7+29$ etc.

\vspace{1cm}

(II) Let $p$ divides $4u_5-1$. Then $4u_5-1=v_5p,~v_5\in
\mathbb{N}$. From (\ref{Erdos-Straus_problem_case3_1st_eqn}), we get
$u_5(y+z)=yzv_5$ i.e. $\frac{v_5}{u_5}=\frac{1}{y}+\frac{1}{z}$. 

(i) If $v_5$ divides $u_5$, then $v_5$ divides $u_5$ and $4u_5-1$. Thus
$v_5=1$, $p=4u_5-1=4u_6+3,~x=(u_6+1)p$ and one possible values of
$y,~z$ are $y=(u_6+2),~z=(u_6+1)(u_6+2)$ where $u_5=u_6+1$.
\begin{equation}
\text{Then}~\frac{4}{p}=\frac{1}{(u_6+1)p}+\frac{1}{(u_6+2)}+\frac{1}{(u_6+1)(u_6+2)}=\frac{1}{\frac{p(p+1)}{4}}+\frac{1}{\frac{p+1}{4}+1}+\frac{1}{\frac{p+1}{4}(\frac{p+1}{4}+1)}. 
\label{Erdos-Straus_problem_p_4n_and_3_3rd}
\end{equation}

(ii) If $v_5=v_6+v_7$ with $v_6$ and $v_7$ divide $u_5$, then
$y=\frac{u_5}{v_6},~z=\frac{u_5}{v_7}$. Also if for some suitable
scalar $r_2\in \mathbb{N}$,
$\frac{1}{y}+\frac{1}{z}=\frac{v_5}{u_5}=\frac{r_2v_5}{r_2u_5}$ with
$r_2v_5=v_8+v_9$ such that $v_8$ and $v_9$ divide $r_2u_5$, then
$y=\frac{r_2u_5}{v_8},~z=\frac{r_2u_5}{v_9}$.

\vspace{1cm}

{\bf Acknowledgment:} {\it I am grateful to the University Grants
  Commission (UGC), New Delhi for awarding the Dr. D. S. Kothari Post
  Doctoral Fellowship from 9th July, 2012 to 8th July, 2015 at Indian
  Institute of Technology (BHU), Varanasi. The self training for this
  investigation was continued during the period.}

\end{document}